\documentclass{article}
\usepackage{PRIMEarxiv}
\usepackage{float}
\usepackage{amsmath}
\usepackage{ upgreek }
\usepackage[utf8]{inputenc} 
\usepackage[T1]{fontenc}    
\usepackage{hyperref}       
\usepackage{url}            
\usepackage{booktabs}       
\usepackage{amsfonts}       
\usepackage{nicefrac}       
\usepackage{microtype}      
\usepackage{lipsum}
\usepackage{fancyhdr}       
\usepackage{graphicx}       
\graphicspath{{media/}}     

\title{Modeling Stock Price Dynamics on the Ghana Stock Exchange: A Geometric Brownian Motion Approach

}

\author{
  D. L. Quayesam \\
  University of Cincinnati \\
  Cincinnati\\
  \texttt{quayesdl@uc.edu} \\
   \And
  A. Lotsi \\
  University of Ghana \\
  Accra\\
  \texttt{alotsi@ug.edu.gh} \\
   \And
    F. O. Mettle \\
  University of Ghana \\
  Accra\\
  \texttt{fomettle@ug.edu.gh} \\
}

\begin{document}
\maketitle

\begin{abstract}

Modeling financial data often relies on assumptions that may prove insufficient or unrealistic in practice. The Geometric Brownian Motion (GBM) model is frequently employed to represent stock price processes. This study investigates whether the behavior of weekly and monthly returns of selected equities listed on the Ghana Stock Exchange conforms to the GBM model. Parameters of the GBM model were estimated for five equities, and forecasts were generated for three months. Evaluation of estimation accuracy was conducted using mean square error (MSE). Results indicate that the expected prices from the modeled equities closely align with actual stock prices observed on the Exchange. Furthermore, while some deviations were observed, the actual prices consistently fell within the estimated confidence intervals. 
\end{abstract}


\section{Introduction}

The evolution of financial asset prices over time is a fundamental aspect of investment, with the ultimate goal being wealth creation. In the world of finance, understanding and predicting stock prices is a fundamental endeavor with far-reaching implications for investors, businesses, and economies at large. Across the globe, the stock market stands as an appealing avenue for investment, and the Ghana Stock Exchange, in particular, garners interest from both local and foreign investors seeking lucrative opportunities. However, achieving expected profits has proven elusive for many investors, often due to challenges in making optimal investment decisions.

Asset prices undergo continuous fluctuations driven by the influx of new information, leading to random changes over short periods. Studies such as \cite{karangwa2008comparing} have observed that both money and commodity markets adhere to the Brownian motion process, reflecting the random nature of market trends. The Geometric Brownian motion is a stochastic process that is widely used in various fields such as finance, physics, and biology. It is characterized by multiplicative noise and can be defined based on different interpretations of stochastic integrals. Among the various approaches to modeling stock prices, the Geometric Brownian Motion (GBM) stands out as a powerful and widely-utilized tool. Rooted in stochastic calculus and probability theory, GBM provides a robust framework for simulating the unpredictable dynamics of financial markets. 

In the face of uncertainties inherent in the equity market, fixed income market, Ghana Alternate market, and general commodity market, individuals seek strategies to make informed decisions and optimize business choices. Given the unpredictable nature of equity prices, investors grapple with the challenge of selecting equities that promise favorable returns.

This study delves into the weekly and monthly behaviors of major equities traded on the Ghana Stock Market, aiming to ascertain whether assumptions proposed by previous researchers regarding the Geometric Brownian Motion (GBM) can be applied to Ghanaian equities. Furthermore, the study endeavors to forecast stock prices and provide valuable insights for investors seeking profitable opportunities on the Exchange. Ultimately, the findings of this study are poised to offer recommendations that could empower investors to capitalize on their knowledge of specific equity paths and potentially enhance profitability in the market.
 
\section{Literature Review} 
Stock price modeling on the Ghana Stock Exchange has evolved significantly over the years, encompassing a diverse range of methodologies and approaches. From traditional econometric techniques to cutting-edge machine learning algorithms, researchers continue to explore innovative ways to understand and predict stock price movements in the dynamic Ghanaian financial landscape. 

Several significant studies have proposed extensions to the Geometric Brownian Motion (GBM) model, introducing jumps and stochastic path volatility processes to enhance its applicability. Merton \cite{merton1976option} pioneered the introduction of the first jump-diffusion model, augmenting the GBM model with a compound Poisson process to represent jumps in stock prices. Dmouj \cite{dmouj2006stock} noted that while GBM appears suitable for ideally modeling future stock prices, practical applications over short time periods often reveal inadequacies. However, Dmouj's study concluded that GBM becomes more accurate when applied over longer time horizons.
In a similar vein, Ladde et al. \cite{ladde2009development} developed a modified stochastic linear model based on classical models to evaluate the precision of the prevailing GBM. Their study subjected stock price data to standard statistical tests to assess GBM's effectiveness. Abidin \cite{abidin2012review} utilized GBM to predict the closing prices of small-sized corporations, suggesting its applicability for estimating daily closing prices up to two weeks ahead.

Mettle  et al. \cite{mettle2014methodology} proposed a methodology for analyzing weekly transaction data from selected equities on the Ghana Stock Exchange (GSE). Their study assumed equity price fluctuations to follow a stochastic process with Markov dependency, categorized into states of price appreciation, stability, and decrement. To facilitate optimal investment decisions on the GSE, the authors proposed decision criteria based on maximum transition probability, minimum mean return time, and maximum limiting distribution of price appreciation for the equities.

 Prempeh et al.\cite{prempeh2022determining} utilized the exponential GARCH model to examine the influence of the COVID-19 pandemic on market volatility. Antwi et al. \cite{antwi2022poisson} developed a stock price model based on pure jump processes, which showed better fit for empirical data compared to diffusion and jump diffusion models. Arhenful et al. \cite{arhenful2021exchange}investigated the link between the Ghanaian exchange rate and stock prices, finding a negative association between the two variables. Akurugu et al. \cite{akurugu2022bayesian} used Bayesian switching volatility models to model the stock returns of GCB bank, identifying two regimes of "turbulent" and "tranquil" market conditions . Gyamfi et al. \cite{gyamfi2021drivers}applied the empirical mode decomposition technique to identify the main drivers of stock prices on the GSE, finding that low frequency and trend components influenced the stock index.

\section{Methodology}

The data used for the study consist of time series of weekly and monthly closing price of 14 randomly selected equities on the GSE. The log returns are used because unlike the simple discrete returns, log returns are time additive \cite{giversen2011continuous} and the statistical properties of the log returns are more tractable \cite{tsay2005analysis}. The return at a particular time $t$ is computed as
\begin{equation}
    R_{t}=ln\frac{p_{t}}{p_{t-1}}
\end{equation}

where $p_{t-1}$ is the stock price at time $(t-1)$ and $p_{t}$  is the stock closing price at time $(t)$.\\
Data exploration begins with plotting log returns against time to assess for trends and unit roots. To complement the visual inspection for stationarity, the Augmented Dickey Fuller (ADF) test is employed to provide empirical evidence. Moreover, normality assumptions of the returns dataset are scrutinized using both graphical and statistical tests. Specifically, histograms, Q-Q plots, and Shapiro-Wilk tests are utilized for this purpose.

Furthermore, this study investigates the potential independence of returns series across different time periods, such as weeks or months. This is accomplished through the application of the Ljung-Box test and the Hurst exponent. The Ljung-Box test assesses the presence of autocorrelation, while the Hurst exponent estimation from empirical data helps evaluate the independence of returns and the existence of randomness within the dataset.

\subsection{Model Specification}
 The stock price process is represented by the stochastic differential equation  
\begin{equation}
dp= \mu p dt+\sigma p ds \;\;\;  0\leq t < T
\end{equation}
where $\mu$ represents the expected rate of return and $\sigma$ denotes volatility of the stock price. If a function $G$ of $p$ follows the Itô’s process, then from Itô’s lemma, we have
\begin{equation}
    dG=\left( \frac{dG}{dp}\mu p + \frac{dG}{dt} + \frac{d^{2}G}{dp^{2}}\right) dt + \frac{dG}{dp}\sigma p ds
\end{equation}
Letting $G=ln p$ and substituting the derivatives of $\frac{dG}{dp}$, $\frac{d^{2}G}{dp^{2}}$ and $\frac{dG}{dt}$ into equation $(3)$, we have

\begin{align*}
    dG & = \left(\mu - \frac{\sigma^{2}}{2}\right) dt + \sigma ds\\
    d(ln p) & = \left(\mu - \frac{\sigma^{2}}{2}\right) dt + \sigma ds
\end{align*}

Integrating the stochastic differential equation above from $t$ to $t+1$, such that $t \leq  t+1 \leq T$

\begin{align*}
ln \frac{p_{t+1}}{p_t} & = \int_{t}^{t+1} \left(\mu - \frac{\sigma^{2}}{2}\right) dt + \int_{t}^{t+1}\sigma ds \\
ln \frac{p_{t+1}}{p_t} & =  \left(\mu - \frac{\sigma^{2}}{2}\right)(t+1-t) + \sigma (s_{(t+1)}-s_{t})
\end{align*}
Taking the expectation of the above equation results
\begin{equation}
E(p_{t+1})=p_{t}\exp{\left(\mu - \frac{\sigma^{2}}{2}\right) \Updelta t}
\end{equation}
where $\Updelta t= (t+1)-t$ and $\Updelta s = s_{t+1}-s_{t}$. We assume that $\Updelta s \sim N(0,\Updelta t)$\\
Computing the second moment and expectation of the stock price in equation $(4)$ above, it can been shown that the variance is
\begin{equation}
    var(p_{t+1})=p_{t}^{2}\exp({2\mu + \sigma^{2}\Updelta t}).(\exp{\sigma{2\Updelta}}-1)
\end{equation}
The $(1-\alpha)100 \%$  confidence interval is constructed using the expected value and the standard deviation of the stock price.
      
\subsection{Parameters Estimation}
 The unknown parameters $\theta =(\mu,\sigma)$ in the process of the GBM are estimated by the method of Maximum Likelihood Estimation (MLE). After the independence and normality assumption have been satisfied, the log likelihood function is stated as
\begin{equation*}
    L(\theta)= \sum_{t=1}^{N} ln\left (f_{\theta}(R_{t})\right)
\end{equation*}

where $f_{\theta}(R_{t})$ represents the probability density function (pdf) of the random variable $R_{t}$ given the parameters $\theta$. 

$f_{\theta}(R_{t}) = \frac{1}{\sigma\sqrt{2\pi(\Updelta t)}}\exp{\left[\frac{-\left[R_{t}-(\mu - \frac{\sigma^{2}}{2})\right]^{2}}{{2\sigma^{2}\Updelta t}} \right]}$

The variable $\sigma$ represents the standard deviation of the log returns and it's often interpreted as the volatility of the asset's returns. The variable $\mu$ represents the drift of the log returns and it's typically the mean rate of return of the asset minus half the variance of the returns.

\subsection{Model Testing}
The data was split as $70\%$ train set and $30\%$ test set. The proposed models for various equities are assessed using the test data set. Expected prices for these equities are calculated, and a comparison is drawn between the forecasted prices from the models and the actual prices observed on the Exchange. The study evaluates the model's performance by computing the Mean Square Error (MSE) and assessing how closely the actual stock prices for different periods align with the estimated 95\% confidence interval. A lower MSE value, typically below 10\% when expressed in percentage terms, indicates a good fit of the data to our model.

\newpage

\section{Results and Discussion}
The methods proposed above are applied to 15 selected equities listed on the Ghana Stock Exchange. This study analysed both their weekly and monthly returns. We specifically choose Ghana Commercial Bank (GCB) equity for demonstration and display its weekly and monthly returns results.

\begin{figure}[ht]
\begin{center}
\includegraphics[width=0.8\columnwidth]{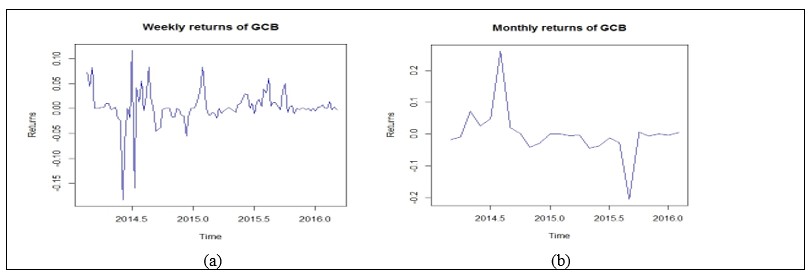}
\end{center}
\caption{Time series plot of weekly and monthly stock returns for GCB}
\label{fig:GCB Trend}
\end{figure}

The GCB equity in figure 1 (a) shows a kind of constant mean and variance since the vertical axis of the graph is stabilized. It is assumed that the weekly return of GCB is stationary. Though the mean and variance in figure 1 (b) looks constant, it is not very clear in the visual check. The Augmented Dickey - Fullers test is shown below 

\begin{table}[ht]
 \caption{Augmented Dickey- Fuller test for GCB}
  \centering
  \begin{tabular}{lll}
    \toprule
    \cmidrule(r){1-2}
    returns     & value     & \textit{p}-value \\
    \midrule
    Weekly & -3.9212  & 0.01566     \\
    Monthly     & -4.2176 & 0.01567      \\
    \bottomrule
  \end{tabular}
  \label{tab:table 1}
\end{table}

The p-value obtained from the Augmented Dickey-Fuller Test for both the weekly and monthly returns of the GCB equity is below the 0.05 significance level. Consequently, we reject the null hypothesis of a unit root (non-stationarity) at the $5\%$ significance level, leading to the conclusion that both the weekly and monthly returns series of GCB are stationary.

Subsequently, the study assessed the normality assumptions using histograms. The histogram graph is depicted in the figure as follows:

\begin{figure}[ht]
\begin{center}
\includegraphics[width=0.8\columnwidth]{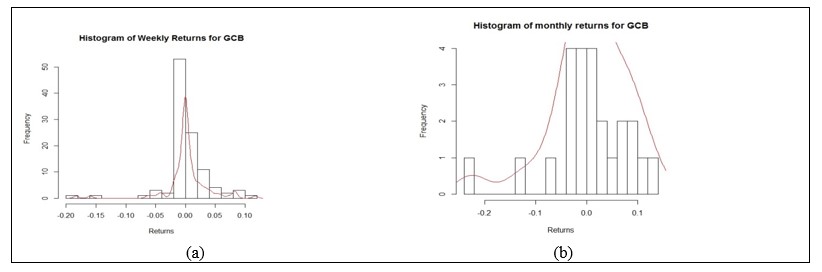}
\end{center}
\caption{Histogram plot of weekly and monthly returns for GCB }
\label{fig:GCB Histogram}
\end{figure}

The histogram representing the weekly returns for GCB in Figure 2(a) does not display a bell shape; it lacks symmetry, and the data points are not evenly distributed around the center. However, the monthly returns of GCB depicted in Figure 2(b) do exhibit a bell shape, with the dataset evenly distributed around the midpoint. To confirm or dispel doubts regarding the normality assumption, the study utilizes a Q-Q plot, as illustrated in the figure below.

\begin{figure}[ht]
\begin{center}
\includegraphics[width=0.8\columnwidth]{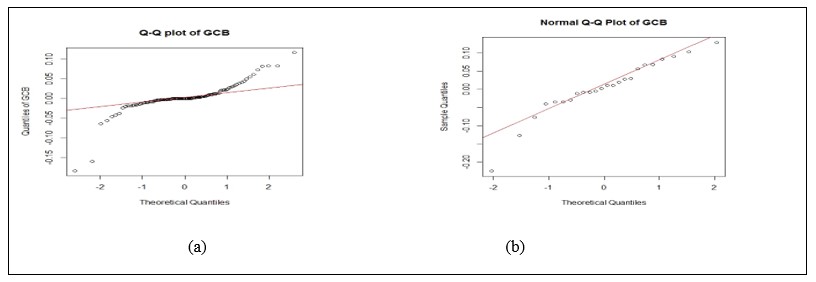}
\end{center}
\caption{Q-Q plot for weekly and monthly returns for GCB}
\label{fig:GCB Q-Q plot}
\end{figure}

The Q-Q plot for the weekly returns of GCB depicted in Figure 3(a) indicates that while the central data points align closely with the line, the tail points deviate from it. This suggests that the weekly returns may not adhere to a normal distribution. Conversely, the Q-Q plot for the monthly returns shown in Figure 3(b) reveals that most points closely follow the straight line, implying a normal distribution.\\
To validate whether the returns conform to a normal distribution, we proceed with the Shapiro-Wilk Test.

\begin{table}[ht]
 \caption{Shapiro-Wilk Test for GCB}
  \centering
  \begin{tabular}{lll}
    \toprule
    \cmidrule(r){1-2}
    returns     & value     & \textit{p}-value \\
    \midrule
    Weekly & 0.76555 & 8.868e-12     \\
    Monthly     & 0.92764 & 0.0863300***     \\
    \bottomrule
  \end{tabular}
  \label{tab:table 2}
\end{table}

In Table 2, the p-value associated with the weekly returns of GCB falls below the chosen significance level of 0.05. Consequently, we reject the null hypothesis, indicating that the weekly returns of GCB equity do not follow a normal distribution. However, the corresponding p-value for GCB monthly returns exceeds the 0.05 significance level, affirming that the monthly returns of GCB adhere to a normal distribution.

To assess the independence within the return series, the study conducted the Ljung-Box test. This test aimed to determine whether the rate of returns in a given week or month correlates with the rate of return from the preceding week or month, respectively.

\begin{table}[ht]
 \caption{Ljung –Box Test for GCB}
  \centering
  \begin{tabular}{lll}
    \toprule
    \cmidrule(r){1-2}
    returns     & value     & \textit{p}-value \\
    \midrule
    Weekly & 1.9975 & 0.1576 ***    \\
    Monthly     & 0.9554 & 0.3284 ***    \\
    \bottomrule
  \end{tabular}
  \label{tab:table 3}
\end{table}

The p-values obtained from Table 3 for the weekly returns (0.1576) and the monthly returns (0.3284) of GCB are both greater than the chosen significance level of 0.05. Consequently, the study does not reject the null hypothesis for either the weekly or monthly returns, suggesting that GCB returns are independent from previous ones.

\newpage

To further explore which equities follow a random process, the study proceeded to estimate the Hurst exponent for both the weekly and monthly returns of GCB.

\begin{table}[ht]
 \caption{Hurst Exponent estimate for GCB}
  \centering
  \begin{tabular}{lll}
    \toprule
    \cmidrule(r){1-2}
    Security     & Weekly     & Monthly \\
    \midrule
     GCB & 0.5538641 & 0.4813197   ***   \\
    \bottomrule
  \end{tabular}
  \label{tab:table 4}
\end{table}

Given that the Hurst exponent for weekly returns exceeds 0.5, the study concludes that GCB's weekly returns do not adhere to a random walk, despite being independent of each other according to the Ljung-Box test. Conversely, examination of Table 4 reveals that the Hurst exponent for monthly returns approaches 0.5. Consequently, based on both the Ljung-Box test for independence and the proximity of the Hurst exponent to 0.5, the study concludes that monthly returns follow a random walk.

From the comprehensive array of tests conducted, it's evident that the monthly returns of GCB conform to all the assumptions of the Geometric Brownian Motion (GBM) model. Conversely, the weekly returns of GCB deviate from the normality and random walk assumptions. Consequently, the same test methods were applied to both weekly and monthly returns of 14 equities listed on the Exchange.

The analysis revealed that the weekly returns of most equities violated at least one assumption of the GBM, with only Tullow Oil (TLW) and Total Oil (TOTAL) adhering to the random walk process. However, among the equities examined, including GCB, five exhibited monthly returns that met all the assumptions of the GBM. The study at this point, focused only on the monthly returns. 

 \subsection{Modeling of Stock prices}
The stochastic model's parameter estimates, as outlined in the methodology, were computed for each monthly return of the five equities that satisfied all the assumptions of the Geometric Brownian Motion (GBM).

 \begin{table}[ht]
 \caption{Parameter estimates }
  \centering
  \begin{tabular}{lll}
    \toprule
    \cmidrule(r){1-2}
    Equity     & $\;\;\;\mu $    &  $\;\;\;\sigma$ \\
    \midrule
    GCB & 0.08499719 & 0.2646442   \\
    GOIL  &-0.26525730 & 0.3090071   \\
    TLW   & 0.06631433 & 0.1683587  \\
    TOTAL  &  0.03200082  & 0.1770301 \\
    UTB   &  0.80071650  & 0.4410409 \\
    \bottomrule
  \end{tabular}
  \label{tab:table 5}
\end{table}
 
The observation of a positive value for the rate of return $(\mu)$ indicates an increase in profit, while a higher volatility $(\sigma)$ suggests rapid fluctuations in the share price. UTB despite having a lower share price compared to the other four equities, exhibits both a higher return and higher volatility. Conversely, GOIL equity is considered risky due to its negative rate of return and high volatility.

The closing price from the last trading date in September 2015 was selected, and a six-month forecast was conducted to assess the predictive capability of our model. By February 2016, despite the relatively broad time period, it was observed that the model's expected price closely approximated the actual price observed on the exchange fell within the confidence interval. 

 \begin{table}[ht]
 \caption{Summary of forecast for Equities}
  \centering
  \begin{tabular}{llll}
    \toprule
    \cmidrule(r){1-2}
    Equity     & Expected Price     &  \;\;\;\;\;\;\;\;\;\;\;CI &  Actual Price \\
               &                    & Lower  \;\;\;\;\;\;\;\;    Upper &         \\
    \midrule
    GCB & 4.20 & 2.65 \;\;\;\;\;\;\;\;\;\; 5.76  & \;\;3.71\\
    GOIL     & 1.66 & 0.94 \;\;\;\;\;\;\;\;\;\; 2.38 &\;\;1.44 \\
    TLW   & 34.00 & 26.04 \;\;\;\;\;\;\;\; 41.96 & \;\;27.93 \\
    TOTAL &   5.49 &  4.14 \;\;\;\;\;\;\;\;\;\; 6.84 & \;\;5.12  \\
    UTB &   0.15 &  0.06 \;\;\;\;\;\;\;\;\;\; 0.24 & \;\;0.12  \\
    
   \bottomrule
  \end{tabular}
  \label{tab:table 6}
\end{table}

The closing price on February 29, 2016, was subsequently utilized as the current price, and forecasts for the expected prices were generated for the final trading days in March, April, and May 2016.

 \begin{table}[ht]
 \caption{Three months forecast of GCB}
  \centering
  \begin{tabular}{llll}
    \toprule
    \cmidrule(r){1-2}
    Month     & Expected Price     &  \;\;\;\;\;\;\;\;\;\;\;CI &  Actual Price \\
               &                    & Lower  \;\;\;\;\;\;\;\;    Upper &         \\
    \midrule
    March & 3.71 & 3.16 \;\;\;\;\;\;\;\;\;\; 4.27  & \;\;3.65\\
    April     & 3.72 & 2.93 \;\;\;\;\;\;\;\;\;\; 4.50 &\;\;3.03 \\
    May   & 3.72 & 2.60 \;\;\;\;\;\;\;\;\;\; 4.84 & \;\;3.02 \\   
   \bottomrule
  \end{tabular}
  \label{tab:table 7}
\end{table}

Table 7 reveals that the monthly expected prices forecasted for GCB remained relatively stable. The forecast for March closely aligns with the actual price, whereas those for April and May exhibit slight deviations. Despite the model's expected prices steadily increase while the actual price of GCB decreases, the observed actual prices for the different months consistently fall within our estimated confidence interval range.

 \begin{table}[ht]
 \caption{Three months forecast of GOIL}
  \centering
  \begin{tabular}{llll}
    \toprule
    \cmidrule(r){1-2}
    Month     & Expected Price     &  \;\;\;\;\;\;\;\;\;\;\;CI &  Actual Price \\
               &                    & Lower  \;\;\;\;\;\;\;\;    Upper &         \\
    \midrule
    March & 1.41 & 1.16 \;\;\;\;\;\;\;\;\;\; 1.72  & \;\;1.45\\
    April     & 1.38 & 1.04 \;\;\;\;\;\;\;\;\;\; 1.72 &\;\;1.51 \\
    May   & 1.35 & 0.94 \;\;\;\;\;\;\;\;\;\; 1.76 & \;\;1.37 \\   
   \bottomrule
  \end{tabular}
  \label{tab:table 8}
\end{table}

According to Table 8, we anticipate a consistent decline in GOIL equity. Despite the forecast closely resembling the actual price, GOIL equity experienced an increase in April before decreasing to 1.37, slightly deviating from our predicted value of 1.35. Nevertheless, the actual prices observed for each month consistently fall within the confidence interval.

 \begin{table}[ht]
 \caption{Three months forecast of TLW}
  \centering
  \begin{tabular}{llll}
    \toprule
    \cmidrule(r){1-2}
    Month     & Expected Price     &  \;\;\;\;\;\;\;\;\;\;\;CI &  Actual Price \\
               &                    & Lower  \;\;\;\;\;\;\;\;    Upper &         \\
    \midrule
    March & 28.08 & 25.41 \;\;\;\;\;\;\;\;\;\; 30.76  & \;\;27.92\\
    April     & 28.24 & 24.43 \;\;\;\;\;\;\;\;\;\; 32.05 &\;\;27.92 \\
    May   & 28.40 & 23.10 \;\;\;\;\;\;\;\;\;\; 33.09 & \;\;27.92 \\   
   \bottomrule
  \end{tabular}
  \label{tab:table 9}
\end{table}

The forecasts outlined in Table 9 for TLW indicate slight increments in the equity. However, our observed actual prices for the various months remained consistent. Both the forecasted and actual prices closely align, with the observed actual prices consistently falling within the confidence range.

 \begin{table}[ht]
 \caption{Three months forecast for TOTAL}
  \centering
  \begin{tabular}{llll}
    \toprule
    \cmidrule(r){1-2}
    Month     & Expected Price     &  \;\;\;\;\;\;\;\;\;\;\;CI &  Actual Price \\
               &                    & Lower  \;\;\;\;\;\;\;\;    Upper &         \\
    \midrule
    March & 5.13 & 4.62 \;\;\;\;\;\;\;\;\;\; 5.65  & \;\;5.10\\
    April     & 5.15 & 4.42 \;\;\;\;\;\;\;\;\;\; 5.88 &\;\;4.90 \\
    May   & 5.16 & 4.26 \;\;\;\;\;\;\;\;\;\; 6.06 & \;\;4.08 \\   
   \bottomrule
  \end{tabular}
  \label{tab:table 10}
\end{table}

The anticipated prices presented in Table 10 indicate a close match between the forecasted closing price for March and the actual price. However, there are slight deviations between the forecasted and actual prices for April and May. Nevertheless, our predicted confidence interval successfully encompasses all the actual prices observed for the various months

 \begin{table}[ht]
 \caption{Three months forecast for UTB}
  \centering
  \begin{tabular}{llll}
    \toprule
    \cmidrule(r){1-2}
    Month     & Expected Price     &  \;\;\;\;\;\;\;\;\;\;\;CI &  Actual Price \\
               &                    & Lower  \;\;\;\;\;\;\;\;    Upper &         \\
    \midrule
    March & 0.13 & 0.10 \;\;\;\;\;\;\;\;\;\; 0.16  & \;\;0.11\\
    April     & 0.14 & 0.09 \;\;\;\;\;\;\;\;\;\; 0.19 &\;\;0.09 \\
    May   & 0.15 & 0.08 \;\;\;\;\;\;\;\;\;\; 0.21 & \;\;0.08 \\   
   \bottomrule
  \end{tabular}
  \label{tab:table 11}
\end{table}

\newpage

Likewise, the anticipated prices shown in Table 11 suggest that the equity price for March closely aligns with the actual price, while there are slight deviations for April and May. Once more, all actual prices for the forecasted months fall within the confidence interval.

\subsection{Mean Square Error of prediction}
The Mean Square Error of prediction revealed that UTB had an error rate of $2.6\%$, GOIL had $6.3\%$, while TLW, GCB, and TOTAL exhibited error rates of $11.9\%$, $32.3\%$, and $41.0\%$ respectively during the three-month model testing. MSE values below $10\%$ indicate the suitability of the model for those equities, with the confidence interval fitness providing justification even for those with higher MSE. 

 \begin{table}[ht]
 \caption{Forecast for monthly returns of the Equities}
  \centering
  \begin{tabular}{llll}
    \toprule
    \cmidrule(r){1-2}
    Equity  & Month   &  Expected Price  &  \;\;\;\;\;\;\;\;\;\;\; CI    \\
               &        &            & \;\; Lower  \;\;\;\;\;\;\;\;    Upper \\
    \midrule
    & July & 3.06 & 2.41 \;\;\;\;\;\;\;\;\;\;\;\; 3.71  \\
     GCB   & August& 3.08 & 2.28 \;\;\;\;\;\;\;\;\;\;\;\; 3.89  \\
        & September& 3.11  & 2.17 \;\;\;\;\;\;\;\;\;\;\;\; 4.04  \\
    \midrule    
    & July & 1.31 & 0.99 \;\;\;\;\;\;\;\;\;\;\;\; 1.64  \\
     GOIL   & August& 1.28 & 0.89 \;\;\;\;\;\;\;\;\;\;\;\; 1.67  \\
        & September& 1.25  & 0.81 \;\;\;\;\;\;\;\;\;\;\;\; 1.70  \\
    \midrule    
    & July & 28.23 & 24.42 \;\;\;\;\;\;\;\;\;\; 32.04  \\
     TLW   & August& 28.39 & 23.69 \;\;\;\;\;\;\;\;\;\; 33.10  \\
        & September& 28.54  & 23.09 \;\;\;\;\;\;\;\;\;\; 34.00  \\
    \midrule    
    & July & 4.10 & 3.52 \;\;\;\;\;\;\;\;\;\;\;\; 4.68  \\
     TOTAL   & August& 4.11 & 3.40 \;\;\;\;\;\;\;\;\;\;\;\; 4.83  \\
        & September& 4.12  & 3.30 \;\;\;\;\;\;\;\;\;\;\;\; 4.95  \\
    \midrule    
    & July & 0.09 & 0.06 \;\;\;\;\;\;\;\;\;\;\;\; 0.12  \\
     UTB   & August& 0.10 & 0.05 \;\;\;\;\;\;\;\;\;\;\;\; 0.14  \\
        & September& 0.10  & 0.05 \;\;\;\;\;\;\;\;\;\;\;\; 0.16  \\
    \bottomrule    
      \end{tabular}
  \label{tab:table 12}
\end{table}

Table 12  provides insights into the expected price movements for the next four months across various equities. GCB equity is anticipated to appreciate in capital gains, with its price not surpassing 4.04. Conversely, GOIL is expected to experience a slight decline, yet its price will not drop below 0.81. For any potential increments, GOIL's price is not expected to exceed 1.70 from its current value. TLW is forecasted to see a modest appreciation, confidently remaining within the range of 23.09 to 34.00. Similarly, the price of TOTAL is projected not to exceed 4.95 from its current value. UTB equity is expected to marginally appreciate by 0.01 at the end of August and remain same in September. 

Among these equities, GCB and TLW present promising options for investors seeking short-term profit within a four-month timeframe.

\section{Conclusions}

The study concludes that investors can leverage the GBM process to inform decision-making and predict stock prices on the Ghana Stock Exchange. Notably, expected stock prices predicted by the model closely align with actual price values, with minor deviations. Nonetheless, all observed actual prices for various months fall within the $95\% $confidence interval. 

The study revealed that among the selected equities, at least one assumption of the Geometric Brownian Motion (GBM) was violated for 13 weekly returns, while the monthly returns of five equities met all GBM assumptions.  This study therefore advocates for the utilization of monthly return prices as a dataset for other researchers, as opposed to daily or weekly returns.This underscores that the behavior of GSE market returns conforms to GBM assumptions for certain equities.

\bibliographystyle{unsrt}  
\bibliography{references}  

\begin{thebibliography}{10}

\bibitem{karangwa2008comparing}
Innocent Karangwa.
\newblock Comparing south african financial markets behaviour to the geometric brownian motion process.
\newblock 2008.

\bibitem{merton1976option}
Robert~C Merton.
\newblock Option pricing when underlying stock returns are discontinuous.
\newblock {\em Journal of financial economics}, 3(1-2):125--144, 1976.

\bibitem{dmouj2006stock}
Abdelmoula Dmouj.
\newblock Stock price modelling: Theory and practice.
\newblock {\em Masters Degree Thesis, Vrije Universiteit}, 2006.

\bibitem{ladde2009development}
GS~Ladde and Ling Wu.
\newblock Development of modified geometric brownian motion models by using stock price data and basic statistics.
\newblock {\em Nonlinear Analysis: Theory, Methods \& Applications}, 71(12):e1203--e1208, 2009.

\bibitem{abidin2012review}
Siti Nazifah~Zainol Abidin and Maheran~Mohd Jaffar.
\newblock A review on geometric brownian motion in forecasting the share prices in bursa malaysia.
\newblock {\em World Applied Sciences Journal}, 17(1):82--93, 2012.

\bibitem{mettle2014methodology}
Felix~Okoe Mettle, Enoch Nii~Boi Quaye, and Ravenhill~Adjetey Laryea.
\newblock A methodology for stochastic analysis of share prices as markov chains with finite states.
\newblock {\em SpringerPlus}, 3:1--11, 2014.

\bibitem{prempeh2022determining}
Kwadwo~Boateng Prempeh, Joseph~Magnus Frimpong, and Newman Amaning.
\newblock Determining the return volatility of the ghana stock exchange before and during the covid-19 pandemic using the exponential garch model.
\newblock {\em SN Business \& Economics}, 3(1):21, 2022.

\bibitem{antwi2022poisson}
Osei Antwi, Kyere Bright, and Martinu Issa.
\newblock Poisson process modeling of pure jump equities on the ghana stock exchange.
\newblock {\em Journal of Applied Mathematics and Physics}, 10(10):3101--3120, 2022.

\bibitem{arhenful2021exchange}
Peter Arhenful, Richard Fosu, and Mathew Owusu-Mensah.
\newblock Exchange rate and stock price nexus: Evidence from ghana.
\newblock {\em Journal of Social and Development Sciences}, 12(4 (S)):9--15, 2021.

\bibitem{akurugu2022bayesian}
Edward Akurugu, Irene~Dekomwine Angbing, Suleman Nasiru, and Abdul~Ghaniyyu Abubakari.
\newblock Bayesian switching volatility models for analysing stock returns in ghana.
\newblock {\em Statistics, Politics and Policy}, 13(2):235--254, 2022.

\bibitem{gyamfi2021drivers}
Emmanuel~N Gyamfi, Frederick~AA Sarpong, and Anokye~M Adam.
\newblock Drivers of stock prices in ghana: an empirical mode decomposition approach.
\newblock {\em Mathematical Problems in Engineering}, 2021:1--7, 2021.

\bibitem{giversen2011continuous}
Jesper Giversen and Mohamed Bendkia.
\newblock Continuous time processes in times of crisis: The case of gbm and cev models.
\newblock 2011.

\bibitem{tsay2005analysis}
Ruey~S Tsay.
\newblock {\em Analysis of financial time series}.
\newblock John wiley \& sons, 2005.

\end{thebibliography}

\end{document}